\documentclass{IEEEtran4PSCC}

\usepackage{xcolor}
\usepackage{booktabs}

\usepackage{cite}

\ifCLASSINFOpdf
   \usepackage[pdftex]{graphicx}
\else
   \usepackage[dvips]{graphicx}
\fi

\usepackage[cmex10]{amsmath}
\usepackage{amsfonts,amssymb,mathtools}
\allowdisplaybreaks

\usepackage[colorlinks,linkcolor=blue,citecolor=red,urlcolor=blue,bookmarks=false,hypertexnames=true]{hyperref} 

\usepackage{algorithm}
\usepackage[noend]{algpseudocode}
 \usepackage{mathabx}

\usepackage{array}

\usepackage{subfigure}

\makeatletter
\let\old@ps@headings\ps@headings
\let\old@ps@IEEEtitlepagestyle\ps@IEEEtitlepagestyle
\def\psccfooter#1{%
    \def\ps@headings{%
        \old@ps@headings%
        \def\@oddfoot{\strut\hfill#1\hfill\strut}%
        \def\@evenfoot{\strut\hfill#1\hfill\strut}%
    }%
    \def\ps@IEEEtitlepagestyle{%
        \old@ps@IEEEtitlepagestyle%
        \def\@oddfoot{\strut\hfill#1\hfill\strut}%
        \def\@evenfoot{\strut\hfill#1\hfill\strut}%
    }%
    \ps@headings%
}
\makeatother

\psccfooter{%
        \parbox{\textwidth}{\hrulefill \\ \small{23rd Power Systems Computation Conference} \hfill \begin{minipage}{0.2\textwidth}\centering \vspace*{4pt} \includegraphics[scale=0.06]{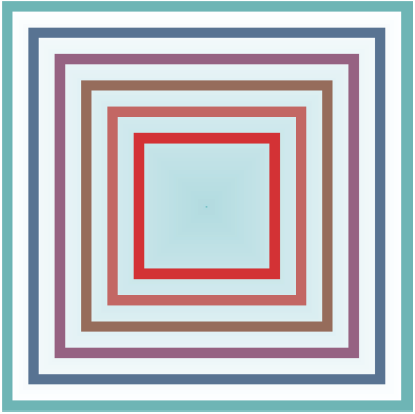}\\\small{PSCC 2024} \end{minipage} \hfill \small{Paris, France --- June 4 -- June 7, 2024}}%
}

\begin{document}
\title{Robust Partitioning and Operation for Maximal Uncertain-Load Delivery in Distribution Grids}

\author{
\IEEEauthorblockN{Hannah Moring$^{1}$, Harsha Nagarajan$^{2}$, Kshitij Girigoudar$^{3}$, David M. Fobes$^{3}$, Johanna L. Mathieu$^{1}$} \\
\IEEEauthorblockA{
$^{1}$ Electrical Engineering and Computer Science, University of Michigan, Ann Arbor, Michigan, USA \\ 
$^{2}$ Applied Mathematics \& Plasma Physics (T-5), Los Alamos National Laboratory, Los Alamos, New Mexico, USA\\
$^{3}$ Information Systems \& Modeling (A-1), Los Alamos National Laboratory, Los Alamos, New Mexico, USA\\
E-mail: \{hmoring, jlmath\}@umich.edu, \{harsha, dfobes\}@lanl.gov, kshitij.girigoudar@gmail.com}
}

\maketitle

\begin{abstract}
To mitigate the vulnerability of distribution grids to severe weather events, some electric utilities use preemptive de-energization as the primary line of defense, causing significant power outages. In such instances, networked microgrids could improve resiliency and maximize load delivery, though the modeling of three-phase unbalanced network physics and computational complexity pose challenges. These challenges are further exacerbated by an increased penetration of uncertain loads. In this paper, we present a two-stage mixed-integer robust optimization problem that configures and operates networked microgrids, and is \textit{guaranteed to be robust and feasible} to all realizations of loads within a specified uncertainty set, while maximizing load delivery. To solve this problem, we propose a cutting-plane algorithm, with convergence guarantees, which approximates a convex recourse function with sub-gradient cuts. Finally, we provide a detailed case study on the IEEE 37-bus test system to demonstrate the economic benefits of networking microgrids to maximize uncertain-load delivery.
\end{abstract}


 \thanksto{This work was supported the U.S. Department of Energy's (DOE) Office of Electricity's Microgrid Research and Development Program under the ``Dynamic Microgrids for Large-Scale DER Integration and Electrification (DynaGrid)'' project and Laboratory Directed Research and Development (LDRD) program under the project ``20230091ER: Learning to Accelerate Global Solutions for Non-convex Optimization''. 
 Moring was also partially supported by NSF Award ECCS-1845093.}

\section{Introduction}
Networked microgrids are gaining significant traction as a means to improve the resilience and economic efficiency of modern distribution grids, especially during extreme weather events and unforeseen contingencies~\cite{Hamidieh2022Microgrids,Li2017Networked,barnes2019resilient}. During such events, isolating affected areas from the rest of the system can be advantageous~\cite{Golari2014Two}. In radial distribution networks, this isolation can result in widespread power outages. Partitioning the network into multiple self-sufficient sub-networks, or ``microgrids'', could alleviate this issue, ideally minimizing load shedding or maximizing load delivery.

Maximizing delivery of \textit{uncertain} loads within distribution networks while minimizing generation costs is challenging. Unlike transmission networks, distribution networks are typically unbalanced. Moreover, if load shedding is necessary, switches are usually employed to shed entire sub-networks rather than just the fraction of load that cannot be met~\cite{fobes2022optimal}. This further complicates balancing power supply and demand.

The integration of distributed energy resources (DERs) has brought  both challenges and opportunities for enhancing  distribution network resiliency~\cite{fobes2022optimal,Zhou2022Three}. Ref.~\cite{Zhou2022Three} presents a Mixed-Integer Linear Program (MILP) with chance constraints to determine a network configuration that minimizes switching costs and expected costs of up-stream power supply while accounting for uncertainty in distributed generation. However,~\cite{Zhou2022Three} does not consider load uncertainty, which is inherent in distribution networks. {\color{black} Some studies do consider load uncertainty. For instance,~\cite{Lee2015Robust} proposes a two-stage robust formulation for distribution network reconfiguration under uncertain loads, while~\cite{Mahdavi2023Robust} presents a single-stage mixed-integer conic program for minimizing losses via distribution network reconfiguration in the presence of load uncertainty. Scenarios wherein load shedding may be necessary are not considered by the above works. With the exception of~\cite{Zhou2022Three}, these papers and many others on network reconfiguration~\cite{babaei2020distributionally,barani2018optimal,arefifar2012supply} only consider balanced networks. Ref.~\cite{Gholami2019Proactive} proposes a two-stage robust formulation for scheduling DERs to mitigate load shedding in the event of distribution network disconnection from the upstream grid, but does not consider network reconfiguration. Although the aforementioned papers model some combination of load uncertainty, load shedding, phase-imbalance, or network reconfiguration, none of them address the unified problem of network reconfiguration aimed at maximizing uncertain-load delivery in unbalanced multi-phase networks.}

To address this planning problem under uncertain loads, we propose a Robust Optimization (RO) approach by modeling loads as uncertainty sets. This approach guarantees a robust partitioning strategy, ensuring the feasibility of the partitioned power network's operation for all load realizations within the specified uncertainty set. 
In this paper, we formulate a two-stage RO model. The first stage is a MILP that optimizes network partitioning decisions while minimizing total unsatisfied loads and generation costs for the nominal load realization. Subsequently, after the load uncertainty is revealed, the second-stage (recourse) problem adjusts generation set-points subject to three-phase, unbalanced power flow constraints. This two-stage problem aims to minimize the worst-case cost of operation and unsatisfied loads, leading to an infinite-dimensional mixed-integer nonlinear optimization problem, typically requiring reformulation for tractability.  

To ensure tractability of the two-stage RO problem, we propose the following solution strategies: {\color{black} (a) employing a flow-based model to enforce radial topology requirements in the first-stage, ensuring all sub-networks are radial, and any islanded sub-network contains an independent voltage/power source, (b) leveraging linearized ``\textsc{LinDist3Flow}'' constraints for tractability in the recourse problem, known for their effectiveness on distribution networks~\cite{robbins2015optimal,Claeys2021No}}, and (c)~generalizing the cutting-plane algorithm from \cite{yang2021robust} to solve the two-stage RO model. We decompose the problem into a mixed-integer master problem and infinitely-many convex subproblems. Given the intractability of handling the latter, we instead find the most violated inequality across all elements of the uncertainty set by solving a $\max-\min$ problem. To optimally solve this $\max-\min$ problem, we exploit the property of convex recourse over a polyhedral set, restricting the search to a finite number of extreme points of the uncertainty set. 

{\color{black} The contributions of this paper are (i) generalization of grid-forming DER constraints from \cite{fobes2022optimal} to ensure any islanded and energized portions of the network contain at least one grid-forming DER,
(ii) development of a tractable, two-stage reformulation of the RO problem to maximize the delivery of uncertain load within an unbalanced distribution network while minimizing generation costs,
(iii) a novel cutting plane algorithm to solve the reformulated RO problem, and
(iv) a detailed numerical study analyzing the sensitivity of planning and operation to uncertainty set parameters, and illustrating the potential benefits of networked microgrids for maximizing load delivery. Additionally, we demonstrate the feasibility of obtained two-stage RO solutions against non-convex AC three-phase power flow through a sampling approach.
} 

\section{Problem Formulation}
In this section, we present necessary notation, the network configuration constraints, and the single-stage robust optimization (RO) problem formulation.

\subsection{Notation and Preliminaries}
\label{subsec:notation}
Bold typeface represents vector notation and blackboard bold typeface represents matrix notation. The $\lvert \cdot \rvert$ operation represents the cardinality when the input is a set, the absolute value when the input is a real number, and magnitude when the input is a complex number.

Consider a distribution network with a set of nodes $\mathcal{N}$, a set of lines $\mathcal{L}$, and a set of phases $\Phi$. There is also a set of transformers $\mathcal{E}^x$ which can be sub-divided into a set of wye-connected transformers $\mathcal{E}^{x,Y}$, and a set of delta-connected transformers $\mathcal{E}^{x,\Delta}$, such that $\mathcal{E}^x=\mathcal{E}^{x,Y}\cup \mathcal{E}^{x,\Delta}$. The network has a meshed structure, but we assume that it can only be operated in radial configurations, as is typical for protection coordination~\cite{Lee2015Robust}. Within the network, there is a set $\mathcal{G}$ of controllable generators (also referred to as DERs) $g$. 
The apparent power injected into the grid by generator $g$ on phase $\phi$ is given by $s^\phi_g$. 

There is also a set of uncontrollable loads $\mathcal{D}$. The apparent power demanded by load $d$ on phase $\phi$ is given by $s^\phi_d = s_d^{0,\phi} + u_{d}^{\phi}$, where $s_d^{0,\phi}$ is the nominal load power and $u_{d}^{\phi}$ represents the load uncertainty defined by the set
\begin{multline}
     \mathcal{U}= \biggl\{ u_{d}^{\phi} \in \mathbb{R}^{|\mathcal{D} \times \Phi|} ~ \bigg| \biggr. ~
      \underline{s}^\phi_{d} - s_{d}^{0,\phi} \leqslant u_{d}^{\phi} \ \leqslant \  \overline{s}^\phi_{d} - s_{d}^{0,\phi}  \\
    ~ \forall \phi \in \Phi, \forall d \in \mathcal{D} \biggr\}.
    \label{eq:uncertaintySet}
\end{multline}
where $\overline{s}^\phi_{d}$ and $\underline{s}^\phi_{d}$ represent the upper and lower bounds of uncertain load $d$ on phase $\phi$.

To capture unbalanced, three-phase power flow typical of distribution grids, the formulation presented in this paper uses the \textsc{LinDist3Flow} model to approximate the AC power flow equations~\cite{Gan2014Convex,Fobes2020PMD}.
Let $v^\phi_i$ represent the voltage magnitude on phase $\phi$ at node $i$ and let $w^\phi_i= (v^\phi_i)^2$. Let $z^{\phi\psi}_{ij} = r^{\phi\psi}_{ij} + jx^{\phi\psi}_{ij}$ represent the self and mutual impedances on the line connecting nodes $i$ and $j$ between phases $\phi$ and $\psi$. The power flowing on the line between nodes $i$ and $j$ on phase $\phi$ is given by $s^\phi_{ij} = p^\phi_{ij} + jq^\phi_{ij}$.

\begin{figure}
    \centering
    \includegraphics[width=0.80\columnwidth]{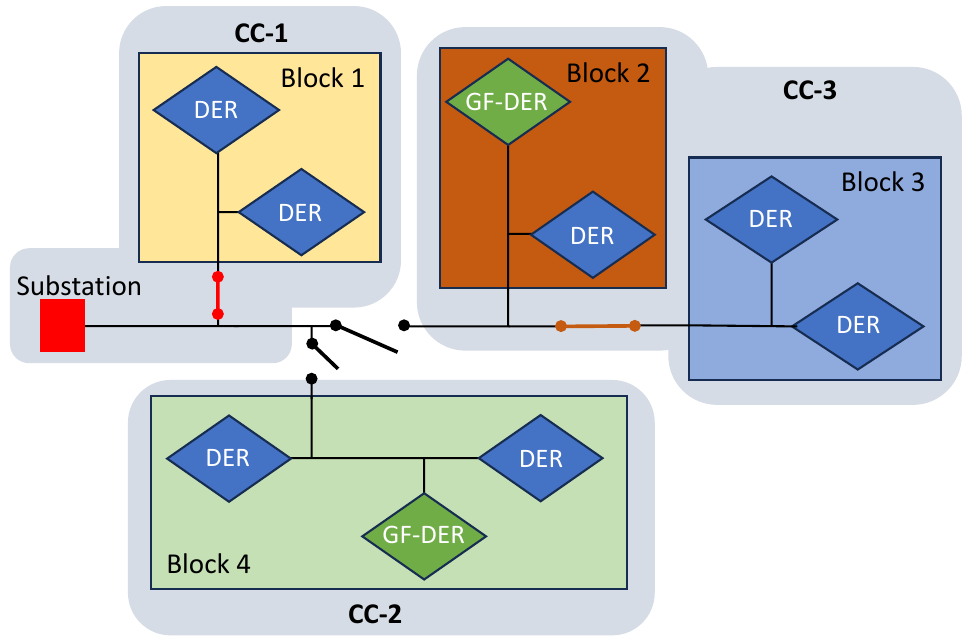}
    \caption{Dynamic partition of networked microgrids using switches. Load blocks are interconnected by switches to form connected components (CCs), acting as independent microgrids. Each CC contains uncertain loads and a set of DERs, with at least one designated as a Grid-Forming DER (GF-DER) when it's CC is energized.}
    \label{fig:load_block}
    \vspace{-.5cm}
\end{figure}

In the network, there is a set of switches $\mathcal{E}^\mathrm{sw}$ that can be opened or closed to reconfigure it. A connected component (CC) is a group of two or more connected nodes~\cite{Lei2020Radiality}.
A microgrid is a CC containing at least one load and at least one DER. Let $\mathcal{B}$ be the set of CCs that exist when every switch in the network is open. The CCs within this set will be referred to as blocks. When a switch between two blocks is closed, the two blocks form one CC. Thus, the number of blocks in the network remains static, while the number of CCs changes with switch configurations, as shown in Fig.~\ref{fig:load_block}.

\subsection{Network Configuration Constraints}
\label{subsec:network_config}
Two general requirements govern network configuration: Firstly, every energized CC must be a spanning tree, devoid of loops. Secondly, each energized CC must include at least one Grid-Forming DER (GF-DER).

\subsubsection{Radiality Constraints}
\label{subsubsec:radial}
We enforce a radial topology in each CC using the directed multi-commodity flow-based model of spanning tree constraints as described in~\cite{Lei2020Radiality}. For brevity, we omit the formulation constraints here and refer the reader to~\cite{Lei2020Radiality} for details. Henceforth, we will refer to this model as the radiality constraints.

\subsubsection{Grid-Forming DER Constraints}
To be energized, a CC must contain a voltage source capable of serving load within it. Distribution network DERs are typically inverter-based resources, which can be either grid-forming or grid-following. A grid-forming DER (GF-DER) controls the AC-side voltage, acting as a voltage source for the network, while a grid-following DER controls the AC-side current and follows the phase angle of the existing grid voltage~\cite{Li2022Revisiting}. {\color{black} In small networks, uncoordinated operation of two or more GF-DERs can lead to independent control of  voltage and frequency, potentially leading to instability~\cite{Sadeque2021Multiple}}. {\color{black} Control strategies facilitating multiple GF-DERs in one microgrid may encounter compatibility challenges with existing older equipment~\cite{Han2016Review}. Thus, the following constraints ensure that each CC contains at least one GF-DER and no more than $k^{\text{der}}$ GF-DER, where $k^{\text{der}} \geq 1$ is a chosen design parameter representing the maximum number of DERs per CC that can simultaneously operate in grid-forming mode. The choice of $k^{\text{der}}$ will depend on equipment capabilities and control strategies.}

We utilize a coloring scheme to determine the presence of a GF-DER 
within a CC based on the implementation in \cite{fobes2022optimal}. Each CC, along with its connected switches, is colored according to the internal block containing the GF-DER. We define binary variables $y^l_{ij}$ to denote whether switch $ij$ is colored by block $l$, $z_l^\mathrm{bl}$ to denote whether block $l$ is energized, and $z_g^\mathrm{inv}$ to denote whether DER $g$ is grid-forming. When switch $ij$ is colored by block $l$, $y^l_{ij}$ equals 1; otherwise, it equals 0. {\color{black} Note that a single switch can be colored by up to $k^{\text{der}}$ blocks.} When a block is energized, $z^\mathrm{bl}$ equals 1; otherwise, it equals 0. Similarly, when a DER is grid-forming, $z^\mathrm{inv}$ equals 1; otherwise, it equals 0.  

Let $\mathcal{G}_l$ and $\mathcal{E}^\mathrm{sw}_l$ be the set of DERs and switches, respectively, connected to block $l$. The relationship between the block states $z_l^\mathrm{bl}$ and the DER states $z_g^\mathrm{inv}$ is defined by
\begin{align}
& |z_i^\mathrm{bl} - z_j^\mathrm{bl}| \leq (1-z^\mathrm{sw}_{ij}), & \forall ij \in \mathcal{E}^\mathrm{sw} \label{eq:ConnectedBlocks} \\
& z^\mathrm{bl}_l \underline{s_g} \leq s_g^\phi \leq z^\mathrm{bl}_l \overline{s_g}, & \forall g \in \mathcal{G}_l, \forall l \in \mathcal{B} \label{eq:BlockGenLimits} \\
& z^\mathrm{bl}_l - \sum_{ij \in \mathcal{E}^\mathrm{sw}_l} z^\mathrm{sw}_{ij} \le \sum_{i \in \mathcal{G}_l} z^\mathrm{inv}_i \le {\color{black}k^{\text{der}}}\mathrlap{z^\mathrm{bl}_l,}  & \forall l \in \mathcal{B} \label{eq:GenPerBlock}
\end{align}
where~\eqref{eq:ConnectedBlocks} states that if switch $ij$ is closed, the blocks it connects are either both energized or both de-energized; \eqref{eq:BlockGenLimits} enforces generation power limits for DERs within an energized block and ensures that DERs within a de-energized block cannot produce power; and \eqref{eq:GenPerBlock} states that, if a block is energized and not connected to another block, then it must contain {\color{black} at least one and no more than $k^{\text{der}}$} GF-DER. 

The coloring scheme used to determine whether a GF-DER exists within a CC is described by
\begin{align}
& \sum_{l \in \mathcal{B}} y^l_{ij} \leqslant {\color{black}k^{\text{der}}} z^\mathrm{sw}_{ij}, \quad  \forall ij \in \mathcal{E}^\mathrm{sw} \label{eq:SwitchColor} \\
\begin{split}
    y^l_{ij} - (1 - z^\mathrm{sw}_{ij}) \le \sum_{i \in \mathcal{G}_l} z^\mathrm{inv}_i \le  {\color{black}k^{\text{der}}} 
    \left(y^l_{ij} + (1 - z^\mathrm{sw}_{ij})\right),\\
   \forall l \in \mathcal{B},\forall ij \in \mathcal{E}^\mathrm{sw}  \label{eq:ColorSwitchGFInv}
\end{split}  \\
\begin{split}
    y^{l'}_{dc} - (1 - z^\mathrm{sw}_{dc}) - (1 - z^\mathrm{sw}_{ab}) \le y^{l'}_{ab} \le  y^{l'}_{dc} + (1 - z^\mathrm{sw}_{dc}) \\ 
    + (1 - z^\mathrm{sw}_{ab}),
    ~\forall ab \ne dc \in \mathcal{E}_l^\mathrm{sw}, \forall l,l' \in \mathcal{B}  \label{eq:ConnectedSwitchesSameColoring} 
\end{split} \\ 
& y_{ij}^l \le \sum_{i \in \mathcal{G}_l} z^\mathrm{inv}_i,~\quad \forall l \in \mathcal{B},\forall ij \in \mathcal{E}^\mathrm{sw} \label{eq:NoInvNoColor} \\
& z^\mathrm{bl}_l \le \sum_{i \in \mathcal{G}_l} z^\mathrm{inv}_i + \sum_{ij \in \mathcal{E}^\mathrm{sw}_l}  \sum_{l \in \mathcal{B}} y^l_{ij},~\quad \forall l \in \mathcal{B} \label{eq:NoColorNoInvNoEnergy} 
\end{align}
where constraint~\eqref{eq:SwitchColor} restricts a closed switch {\color{black} to at most $k^{\text{der}}$ colors}, while an open switch has no color. Constraint~\eqref{eq:ColorSwitchGFInv} ensures that if switch $ij$ is closed and colored by block $l$, then it is in the same CC as block $l$, which contains {\color{black} at least one} GF-DER. Constraint~\eqref{eq:ConnectedSwitchesSameColoring} ensures that if two switches are closed and connected to the same block, they must have the same color. Constraint~\eqref{eq:NoInvNoColor} enforces that a switch cannot have color $l$ unless block $l$ contains {\color{black} at least one} GF-DER. Constraint~\eqref{eq:NoColorNoInvNoEnergy} states that a block is de-energized unless it contains {\color{black} at least one} GF-DER or it is connected to {\color{black} at least one} closed switch.

Let $\mathcal{E}^v_l$ be a set of virtual edges between block $l$ and all other blocks $l'$, oriented from $l$ to $l'$ (there are $|\mathcal{B}|-1 $ of these edges for each $l$). Let $\xi^l_{ij}$ represent a unit of flow of commodity $l$ across the virtual edge from $i$ to $j$. For each $ij \in \mathcal{E}^\mathrm{sw}_l$, there is an arbitrary orientation (e.g., either $i=l$ or $j=l$). Then the following set of constraints can be used to determine if there is a connection between a GF-DER and a CC.
\begin{align}
\begin{split}
    -z^\mathrm{sw}_{ij} (|\mathcal{E}^\mathrm{sw}|) \le \eta_{ij}^l\le z^\mathrm{sw}_{ij} (|\mathcal{E}^\mathrm{sw}|), \ \ 
   \forall ij \in {\cal E}_l^\mathrm{sw}, \ \forall l \in \mathcal{B} \label{eq:RealLineFlow}
\end{split} \\
& 0 \le \xi_{ij}^l \le 1 ,~\quad \forall ij \in \mathcal{E}^{v}_l, \forall l \in \mathcal{B} \label{eq:VirtualLineFlow} \\
\begin{split}
    \sum_{ij \in \mathcal{E}^\mathrm{sw}_l : i = l} \eta_{ij}^l - \sum_{ij \in \mathcal{E}^\mathrm{sw}_l : j = l} \eta_{ij}^l + \sum_{ij \in \mathcal{E}^{v}_l} \xi_{ij}^l = |\mathcal{B}| - 1 ,\\
    \forall l \in \mathcal{B} \label{eq:FlowForEveryBlock}
\end{split}  \\
\begin{split}
    \sum_{ij \in \mathcal{E}^\mathrm{sw}_l : i = l'} \eta_{ij}^l - \sum_{ij \in \mathcal{E}^\mathrm{sw}_l : j = l'} \eta_{ij}^l  -  \xi_{ll'}^l = -1, \\
    \forall l' \ne l, \forall l \in \mathcal{B} \label{eq:OneFlowConsumed}
\end{split}  \\
& y_{ij}^l \le 1 - \xi_{ll'}^l ,~\quad \forall l' \ne l, ij \in \mathcal{E}^\mathrm{sw}_{l'}, \forall l \in \mathcal{B} \label{eq:VirtualFlowNoColor} \\
& z^\mathrm{bl},z^\mathrm{sw},z^\mathrm{inv},y^{l} \in \{0,1\} \label{eq:binaries}
\end{align}
Constraint~\eqref{eq:RealLineFlow} limits the flow magnitude $\eta_{ij}^l$ across a closed switch $ij$ from block $l$ by the number of switches in the network. 
Constraint~\eqref{eq:VirtualLineFlow} restricts flow across each virtual edge to at most one unit. Constraints~\eqref{eq:FlowForEveryBlock} and \eqref{eq:OneFlowConsumed} collectively ensure flow balance among blocks $l$ and $l'$.
Constraint~\eqref{eq:VirtualFlowNoColor} ensures $y^l_{ij}=0$ for every $ij \in \mathcal{E}^\mathrm{sw}_{l'}$ if there is flow on the virtual line between $l$ and $l'$. Lastly~\eqref{eq:binaries} defines binary variables $z^\mathrm{bl},z^\mathrm{sw}, z^\mathrm{inv}$ and $y^l$.

\subsection{Single-Stage Robust Partitioning and Operation Problem}
\label{subsec:single_stage}
We now define the single-stage Robust Partitioning and Operation Problem (RPOP) for distribution grids. Its optimal solution comprises of network partitions and generator/DER set-points that minimize load shed and operational costs while ensuring feasibility and robustness against all realizations of uncertain loads. The formulation is as follows:
\begin{subequations} 
\label{eq:FullForm}
\begin{align}
     \min~ & \sum_{l\in\mathcal{B}}  \alpha_l(1-z^\mathrm{bl}_l)+ \sum_{\substack{g \in G}} \sum_{\substack{\phi \in \Phi_g}}c_{1,g}s^\phi_{g} + c_{0,g} \label{obj} \\
     \mathrm{s.t.}~ & \text{Radiality~Constraints} ~(\text{Sec.}~ \ref{subsubsec:radial}),~\eqref{eq:ConnectedBlocks}-\eqref{eq:binaries} \\ 
    \begin{split}
         z^\mathrm{bl}_l\underline{v_{i}}^2 \leq w_{i,\phi} \leq z^\mathrm{bl}_l\overline{v_{i}}^2, ~\forall \phi \in \Phi, \\
         \forall i \in \mathcal{N}_l, \forall l \in \mathcal{B} \label{eq:voltLim}
     \end{split}  \\
     \begin{split}
         \mathbf{w}_i = \mathbf{w}_j - \mathbb{M}^p_{ij}\mathbf{p}_{ij} - \mathbb{M}^q_{ij}\mathbf{q}_{ij}, \\
         \forall i,j \in \mathcal{N}, ~\forall ij \in \mathcal{L} \label{eq:VoltBal}
     \end{split}  \\
     & \lvert s^\phi_{ij} \rvert \leq \overline{s_{ij}^{\phi}}, ~\quad \forall \phi \in \Phi_{ij}, \forall ij \in \mathcal{L} \cup \mathcal{E}^{x} \label{eq:PowerLim}\\
     & \lvert s^\phi_{ij} \rvert \leq z^\mathrm{sw}_{ij}\overline{s_{ij}^{\phi}}, ~\quad \forall \phi \in \Phi_{ij}, \forall ij \in \mathcal{E}^\mathrm{sw} \label{eq:SWPowerLim}\\
     & w^\phi_i = n_{ij}^2 w^\phi_j, ~\quad \forall \phi \in \Phi_{ij}, \forall ij \in \mathcal{E}^{x,Y} \label{eq:XfmrYV}\\
     & s^\phi_{ij} = s^\phi_{ji}, ~\quad \forall \phi \in \Phi_{ij},\forall ij \in \mathcal{E}^{x,Y} \label{eq:XfmrY} \\
     \begin{split}
         & 3(w_{i,\phi}+w_{i,\psi})= 2(n_{ij})^2w_{j,\phi}, \\
         &\quad \quad \forall (\phi,\psi) \in \{(a,b),(b,c),(c,a)\}, ~\forall ij \in \mathcal{E}^{x,\Delta} \label{eq:XfmrDeltaV}
     \end{split}\\
    \begin{split}
        2p_{ij,\phi}=-(p_{ji,\phi}+p_{ji,\psi})+(q_{ji,\psi}-q_{ji,\phi})/\sqrt{3}, \\ ~\forall (\phi,\psi) \in \{(a,c),(b,a),(c,b)\}, ~\forall ij \in \mathcal{E}^{x,\Delta} \label{eq:XfmrDeltaP}
    \end{split}\\
    \begin{split}
        2q_{ij,\phi}=(p_{ji,\phi}-p_{ji,\psi})/\sqrt{3}-(q_{ji,\psi}+q_{ji,\phi}), \\
        ~\forall (\phi,\psi) \in \{(a,c),(b,a),(c,b)\}, ~\forall ij \in \mathcal{E}^{x,\Delta} \label{eq:XfmrDeltaQ}
    \end{split} \\
    \begin{split}
         \sum_{ij \in \mathcal{E}_i} s^\phi_{ij} = \sum_{g \in i} s^\phi_{g}
         - z^\mathrm{bl}_k \sum_{d \in \mathcal{D}_i}s^\phi_d 
         - \sum_{c \in \Phi_i} y^\phi_c w^\phi_i, \\ ~\forall \phi \in \Phi_i,\forall i \in l,~\forall l \in \mathcal{B} \label{eq:PowerBalance}
     \end{split} 
\end{align}
\end{subequations}
where $\mathcal{E} = \mathcal{L}\cup\mathcal{E}^\mathrm{sw}\cup\mathcal{E}^{x}$ is the set of all lines, switches, and transformers.
The objective function minimizes the number of de-energized blocks and the generation cost, where $\alpha_l$ is a weighting parameter for block $l$, indicating block priority, and $c_{1,g}, c_{0,g}$ are parameters of the linear cost function of generator $g$. The voltage limits at each bus are enforced by~\eqref{eq:voltLim} and the linearized voltage drop between bus $j$ and downstream bus $i$ is computed in~\eqref{eq:VoltBal}, also known as the \textsc{LinDist3Flow}~\cite{Gan2014Convex,robbins2015optimal}, where the variables represented with bold typeface are 3-phase vectors, and 
\begin{align*}
    \mathbb{M}^p_{ij} &= \begin{bmatrix}
    -2r_{aa} & r_{ab}-\sqrt{3}x_{ab} & r_{ac}+\sqrt{3}x_{ac}\\ r_{ba}+\sqrt{3}x_{ba} & {-}2r_{bb} & r_{bc}-\sqrt{3}x_{bc}\\ r_{ca}-\sqrt{3}x_{ca} & r_{cb}+\sqrt{3}x_{cb} & {-}2r_{cc} \end{bmatrix}_{ij}, \\
    \mathbb{M}_{ij}^{q} &=\begin{bmatrix} {-}2x_{aa} & x_{ab}+\sqrt{3}r_{ab} & x_{ac}-\sqrt{3}r_{ac}\\ x_{ba}-\sqrt{3}r_{ba} & {-}2x_{bb} & x_{bc}+\sqrt{3}r_{bc}\\ x_{ca}+\sqrt{3}r_{ca} & x_{cb}-\sqrt{3}r_{cb} & {-}2x_{cc} \end{bmatrix}_{ij}.
\end{align*}
Constraints~\eqref{eq:PowerLim} and ~\eqref{eq:SWPowerLim}, 
define the power flow limits for each line, transformer, and switch. 
Constraints~\eqref{eq:XfmrYV} and~\eqref{eq:XfmrDeltaV} define the voltage transformations across each wye-connected and delta-connected transformer, where $n_{ij}$ is the tap ratio of transformer $ij$. Constraints~\eqref{eq:XfmrY},~\eqref{eq:XfmrDeltaP} and~\eqref{eq:XfmrDeltaQ} define the relationship between the directions of power flow on each wye-connected and delta-connected transformer. Finally,~\eqref{eq:PowerBalance} is the linearized power balance equation.

\section{A Cutting-Plane Algorithm to solve the RPOP}
Since it is intractable to solve the RPOP as given in Section~\ref{subsec:single_stage} using any off-the-shelf solvers, we propose a two-stage version of the RPOP, which can be solved effectively to global optimality using an iterative cutting-plane algorithm. 

\subsection{Two-stage Robust Partitioning and Operation Problem}
\label{subsec:2stage_RPOP}
The RPOP in \eqref{eq:FullForm} can be exactly reformulated as a master problem and a set of subproblems~\cite{Geoffrion1972Generalized}. The master, or first-stage, decision variables include controllable power injections, $\mathbf{s}^\phi_g~ \forall \phi \in \Phi,~\forall g \in \mathcal{G}$, switch configurations $z^\mathrm{sw}_{ij},~ \forall ij \in \mathcal{E}^\mathrm{sw}$,  generator operating states $z^\mathrm{inv}_{g},~ \forall g \in \mathcal{G}$ as either grid-forming or grid-following, block energized states $z^\mathrm{bl}_l,~ \forall l\in \mathcal{B}$, and additional variables related to network configuration constraints.
Let $\mathbf{x}$ represent a vector of these variables $\mathbf{s}^\phi_g,z^\mathrm{sw}_{ij},z^\mathrm{inv}_{g},z^\mathrm{bl}_l$ including every generator, switch, and block in the network. Let $\mathbf{s}_{d}^*$ be the worst-case uncertain load realization.
The master problem, which relaxes the single-stage problem by ignoring the power flow constraints, is
\begin{subequations} 
\begin{align}
    \min~ & \sum_{l\in\mathcal{B}}  \alpha_l(1-z^\mathrm{bl}_l)+ \sum_{\substack{g \in G}} \sum_{\substack{\phi \in \Phi_g}}c_{1,g}s^\phi_{g} + c_{0,g} + \theta \nonumber \tag{$M$}\label{eq:master}\\
    \mathrm{s.t.}~ & \mathrm{Radiality~Constraints}~(\text{Sec.}~ \ref{subsubsec:radial}),~\eqref{eq:ConnectedBlocks}-\eqref{eq:binaries} \\ 
    & V_2(\mathbf{x}^*_k,\mathbf{s}_{d,k}^*) + \mathbf{\pi}_k^{\intercal} A(\mathbf{x}-\mathbf{x}^*_k) \leq \theta,
        \forall k=1,2,... \label{eq:Cuts} 
\end{align} 
\end{subequations} 
where a subscript $k$ denotes the value at the $k^\mathrm{th}$ iteration, $\mathbf{x}^*_k$ denotes the optimal solution of~\eqref{eq:master} at the $k^\mathrm{th}$ iteration, $V_2(\mathbf{x}^*_k,\mathbf{s}_{d,k}^*)$ is the objective value of the $k^\mathrm{th}$ subproblem, $\mathbf{\pi}$ is a vector of dual variables corresponding to equations in~\eqref{eq:subproblem} containing $\mathbf{x}$, and $A$ is the coefficient on $\mathbf{x}$ in those equations.

In the second stage of the problem (or the subproblem), adjustments to the controllable power injections are permitted after the uncertainty is revealed. These adjustments are represented by variables $o^{+}$ and $o^{-}$, bounded by generator ramping capabilities and capacity constraints. Additionally, squared bus voltages $w^\phi_i~ \forall \phi \in \Phi_i,\forall i \in \mathcal{N}$, and power flows $s^\phi_{ij}~ \forall \phi \in \Phi_{ij},\forall ij \in \mathcal{E}$ are determined by the power flow equations in this stage. The subproblem formulation is
\begin{subequations} 
\label{eq:subprob}
\begin{align}
    &\min~ \omega \left( \sum_{i \in \mathcal{N}}\sum_{\phi \in \Phi_i}\left(h^{+,\phi}_i + h^{-,\phi}_i \right) \right)
    + \sum_{g \in G}\sum_{\phi \in \Phi_g} c_{1,g} o^\phi_{g} \nonumber \tag{$S_1$}\label{eq:subproblem}\\
    & s.t.~ \quad \eqref{eq:voltLim}-\eqref{eq:XfmrDeltaQ}, \nonumber \\
    \begin{split}
         \sum_{ij \in \mathcal{E}_i} s^\phi_{ij} + h^{+,\phi}_i -  h^{-,\phi}_i = \sum_{g \in i} \left(s^{\phi*}_{g} + o^{+,\phi}_g - o^{-,\phi}_g \right)
         \\ - z^\mathrm{bl*}_k \sum_{d \in \mathcal{D}_i}s^\phi_d 
         - \sum_{c \in \Phi_i} y^\phi_c w^\phi_i, ~\forall i \in k,~\forall k \in \mathcal{B}, \label{eq:SubPowerBalance}
     \end{split} \\
     & o^{+,\phi}_{g} \leq z^\mathrm{bl*}_i\overline{s^\phi_{g}}  - s^{\phi*}_{g}, ~ \forall \phi \in \Phi_g, \forall g \in i, \forall i \in \mathcal{B},  \label{eq:SubModUpLims}\\
     & o^{-,\phi}_{g} \leq s_{g}^{\phi*} - z^\mathrm{bl*}_i\underline{s^\phi_{g}}, ~ \forall \phi \in \Phi_g,~\forall g \in i,~ \forall i \in \mathcal{B}, \label{eq:SubModDownLims}\\
     & o^{+,\phi}_{g} \leq \overline{o^\phi_{g}}, ~ \forall \phi \in \Phi_g,~\forall g \in i,~ \forall i \in \mathcal{B}, \label{eq:Sub+RampLims}  \\
     & o^{-,\phi}_{g} \leq \overline{o^\phi_{g}}, ~ \forall \phi \in \Phi_g,~\forall g \in i,~ \forall i \in \mathcal{B}, \label{eq:Sub-RampLims} \\
     & o^{+,\phi}_{g}, o^{-,\phi}_{g}\geq 0, ~ \forall \phi \in \Phi_g,~\forall g \in i,~ \forall i \in \mathcal{B}, \label{eq:SubModPos} \\
     & h^{+,\phi}_i, h^{-,\phi}_i, \geq 0 , ~ \forall \phi \in \Phi_i,~ \forall i \in \mathcal{N}. \label{eq:SubSlackPos}
\end{align}
\end{subequations}
where $\omega$ is a weighting parameter, and $h_i^{+,\phi}, h_i^{-,\phi}$ are slack variables indicating a power balance violation on phase $\phi$ at node $i$. If any of these slack variables are non-zero at optimality, there are no feasible $o^{+},o^{-}$ that satisfy the power flow equations for the candidate solution $\mathbf{x}^*$ and uncertain load realization $\mathbf{s}_d$. In such cases, a feasibility cut (\ref{eq:Cuts}) is generated to ensure that $\mathbf{x}^*$ 
cannot be chosen again in subsequent iterations of the master problem.
Constraints~\eqref{eq:SubModUpLims}-\eqref{eq:SubModPos} set limits on generator set-point modifications based on ramping and capacity limits.

The decomposition algorithm cannot be directly implemented due to infinitely many subproblems arising from uncertain load realizations. Instead, it is preferable to identify and solve only the subproblem with most violated power balance constraints within the given uncertainty set. This entails addressing a bilevel $\max-\min$ problem of the form:
\begin{subequations}  
\begin{align}
    & \max_{\mathcal{U}} \ \min ~ \omega \left( \sum_{i \in \mathcal{N}}\sum_{\phi \in \Phi_i}\big(h^{+,\phi}_i + h^{-,\phi}_i \big) \right) + \sum_{g \in G}\sum_{\phi \in \Phi_g} c_{1,g} o^\phi_{g} \nonumber \tag{$S_2$}\label{eq:MaxMin} \\ 
    & s.t. \quad \eqref{eq:voltLim}-\eqref{eq:XfmrDeltaQ}, \eqref{eq:SubModUpLims}-\eqref{eq:SubSlackPos} \nonumber
\end{align}
\end{subequations} 
According to convex optimization theory, the optimal solution to the $\max-\min$ problem~\eqref{eq:MaxMin} will manifest at one of the extreme points of the uncertainty set since the inner minimization problem constitutes a convex function of the uncertain load parameters~\cite{Bertsekas2009Convex}. Leveraging this result, we can substitute $\mathcal{U}$ in~\eqref{eq:MaxMin} with a finite-dimensional uncertainty set 
\begin{multline*}
     \widehat{\mathcal{U}}= \biggl\{ u_{d}^{\phi} \in \mathbb{R}^{|\mathcal{D} \times \Phi|} ~ \bigg| \biggr. ~
      u_{d}^{\phi} = \zeta_d^{+}(\overline{s}^\phi_{d}) + \zeta_d^{-}(\underline{s}^\phi_{d}) - s_{d}^{0,\phi},      \\
      \hspace{0.8cm} \zeta_d^{+}+\zeta_d^{-}=1, \ \zeta_d^{+},\zeta_d^{-} \in \{0,1\}, 
    ~ \forall \phi \in \Phi, \forall d \in \mathcal{D} \biggr\}.
\end{multline*}
To simplify the uncertainty modeling further, we assume that the load variation for multi-phase loads at a bus is uni-directional, meaning $\zeta_d^{+},\zeta_d^{-}$ are indexed by loads but not by phases.
Thus, the reduced finite-dimensional~\eqref{eq:MaxMin} problem can be replaced by solving the modified subproblem for every possible combination of $\zeta_d^+,\zeta_d^-$, selecting the solution that yields the largest objective, i.e., 
\begin{subequations} 
\label{eq:ExtremeSub}
\begin{align}
    \begin{split}
        \min~ \omega \left( \sum_{i \in \mathcal{N}}\sum_{\phi \in \Phi_i}\left(h^{+,\phi}_i + h^{-,\phi}_i \right) \right)
    + \sum_{g \in G}\sum_{\phi \in \Phi_g} c_{1,g} o^\phi_{g}
    \end{split} \tag{$E$}\label{eq:enumerate} \\
    & s.t. \quad \eqref{eq:voltLim}-\eqref{eq:XfmrDeltaQ}, \eqref{eq:SubModUpLims}-\eqref{eq:SubSlackPos}, \nonumber \\
    \begin{split}
         \sum_{ij \in \mathcal{E}_i} s^\phi_{ij}  + h^{+,\phi}_i -  h^{-,\phi}_i = \sum_{g \in i} \left(s^{\phi*}_{g} + o^{+,\phi}_g - o^{-,\phi}_g \right)
         \\ - z^\mathrm{bl*}_k \sum_{d \in \mathcal{D}_i} \big( s^{0,\phi}_d 
          + \zeta^+_d\overline{s}^{\phi}_d - \zeta^-_d\underline{s}^{\phi}_d \big) - \sum_{c \in \Phi_i} y^\phi_c w^\phi_i, \\ ~\forall i \in k,~\forall k \in \mathcal{B}, \nonumber
     \end{split} \\
     & \zeta^+_d + \zeta^-_d  = z^\mathrm{bl*}_i, \ \ \zeta^+_d,\zeta^-_d \in \{0,1\}, \ \forall d \in \mathcal{D}_i,~ \forall i \in\mathcal{B}. \nonumber
\end{align}
\end{subequations}
This approach still requires solving $2^{|\mathcal{D}|}$  number of quadratically-constrained convex programs \eqref{eq:enumerate}, and may not scale efficiently as the number of uncertain loads increases.  

\subsection{A Cutting-plane Algorithm}
We propose Algorithm~\ref{alg:cutting_plane} to solve the two-stage RPOP. In summary, the master problem~\eqref{eq:master} is solved for the network configuration and generator set-points (master solution, $\mathbf{x}^*$) in Step~\ref{alg:solve_M}. Given this $\mathbf{x}^*$, \eqref{eq:MaxMin} is solved for every scenario in $\widehat{\mathcal{U}}$ to get the worst-case load realization in Step~\ref{alg:solve_maxmin}. Finally, the subproblem~\eqref{eq:subproblem} is solved in Step~\ref{alg:solve_S} to check if the chosen master solution is feasible for the worst-case load realization. If they are feasible, the algorithm terminates with an optimal solution to the RPOP. Otherwise, if there is non-zero slack $h^{+,\phi},  h^{-,\phi}$, then a sub-gradient cut is added to~\eqref{eq:master} in Step~\ref{alg:add_cut} and the process is repeated.
\noindent
\\
\textbf{Convergence guarantee:} 
Algorithm~\ref{alg:cutting_plane} ensures convergence within a finite number of iterations due to the finite number of load blocks in any distribution grid. In the worst-case scenario, the sub-gradient cuts will necessitate de-energizing every block in the network to meet the termination criterion, as indicated in step~\eqref{alg:termination} of Algorithm \ref{alg:cutting_plane}.

\begin{algorithm}
  \caption{: Cutting-plane algorithm for two-stage RPOP}
  \label{alg:cutting_plane}
\begin{algorithmic}[1]
\State Initialize: $k \leftarrow 1$, $h^{+,\phi}_i \leftarrow \infty$, $h^{-,\phi}_i \leftarrow \infty~\forall i \in \mathcal{N}$, $\epsilon > 0$, $V_2(\mathbf{x}_1^*,\mathbf{s}_{d,1}^*)\leftarrow 0$, and $\mathbf{\pi}_1 \leftarrow 0$ 
\While {$\sum_{i \in \mathcal{N}}\sum_{\phi \in \Phi_i}\left(h^{+,\phi}_i + h^{-,\phi}_i \right) > \epsilon$} \label{alg:termination}
\State Append $V_2(\mathbf{x}_k^*,\mathbf{s}_{d,k}^*) + \mathbf{\pi}_k^{\intercal} A(\mathbf{x}-\mathbf{x}_k^*) \leq \theta$ to~\eqref{eq:master} \label{alg:add_cut}
\State Solve master problem~\eqref{eq:master} for $\mathbf{x}_k^*$ \label{op1} \label{alg:solve_M}
\State  Given $\mathbf{x}_k^*$, find the worst-case uncertainty realization $s_{d,k}^*$ by solving~\eqref{eq:enumerate} for every scenario in $\hat{\mathcal{U}}$ \label{alg:solve_maxmin}
\State Given $\mathbf{x}_k^*$ and $s_{d,k}^*$, solve Subproblem~\eqref{eq:subproblem} for $h^{+,\phi},  h^{-,\phi}, o^{+,\phi}, o^{-,\phi}$ \label{alg:solve_S}
\State Set $k \leftarrow k + 1$
\EndWhile
\State Output: $\mathbf{x}_{k-1}^*, o^{+,\phi}, o^{-,\phi}$
\end{algorithmic}
\end{algorithm}

\section{Numerical Results}
\label{sec:results}
We now present results illustrating how robust partitioning and operation of networked microgrids using Algorithm~\ref{alg:cutting_plane} can maximize load delivery.

\subsection{Implementation details}
The proposed optimization formulations and Algorithm~\ref{alg:cutting_plane} were implemented using JuMP v1.13.0~\cite{dunning2017jump} in the Julia v1.9.2 programming language. PowerModelsONM.jl~\cite{fobes2022optimal} served as the foundational framework for all implementations. All optimization tasks were executed with Gurobi 10.0.3~\cite{gurobi} as an MILP solver on hardware comprising a 1.3 GHz 4-Core Intel Core i7 processor with 16GB memory.

\subsection{Test case: IEEE 37-bus test system}
\label{subsec:test_case}
\begin{figure}
\centering
\includegraphics[width=3.3in]{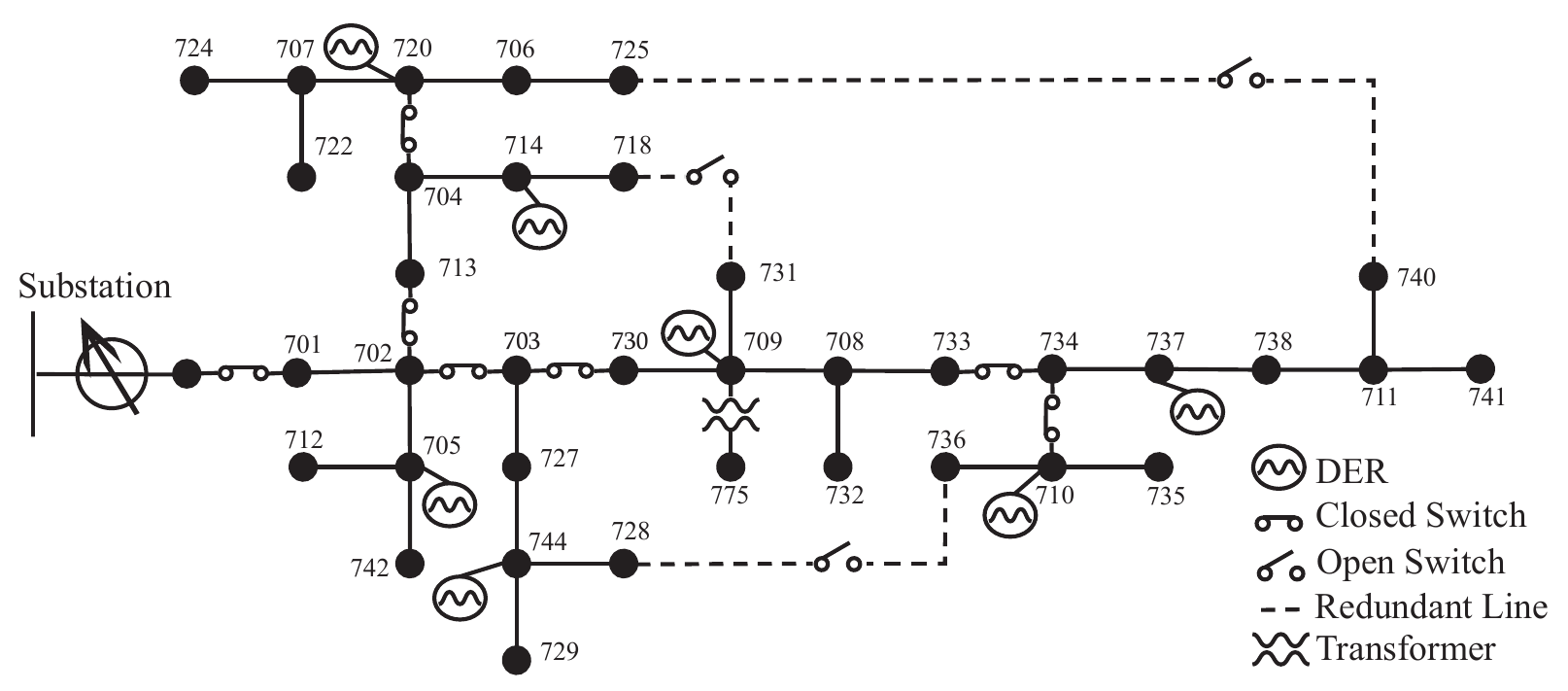}
\vspace{-2mm}
\caption{Single-line diagram of the 37-bus network without loads shown.}
\label{fig:37bus_network}
\vspace*{-2mm}
\end{figure}

We demonstrate the benefits of robust network partitioning against load uncertainty using a modified version of the unbalanced, three-phase, \textit{IEEE 37-bus test network} \cite{kersting2001radial}, whose single-line diagram is shown in Fig.~\ref{fig:37bus_network}.  This modified network includes 10 controllable switches and 7 controllable  DERs. Among these switches, three are part of redundant lines (dashed) with sufficient capacity to facilitate power transfer between networked microgrids. These controllable switches allow the network to be partitioned into a maximum of 7 blocks. We consider the case where a disruptive event occurs upstream of the network's substation, necessitating isolation from the main grid. To achieve this isolation, the switch between the substation and node 701 is opened. 

The network comprises 30 loads, distributed non-uniformly across nodes and three phases, resulting in significant phase imbalance~\cite{kersting2001radial}. The total nominal load in the network is 2542~kW. Among these loads, six are subject to modeling uncertainty, distributed across the network, while the remaining 24 loads are held fixed at their nominal values. We consider using this reduced set of uncertain loads reasonable, owing to the computational burden of solving the exponential number of~\eqref{eq:enumerate} problems for each iteration of the algorithm. 

{\color{black} We assume no autonomous coordination among DERs, necessitating that each energized CC in the network has precisely one GF-DER ($k^{der} =1$), with all other DERs in the CC designated as grid-following.} The total generation capacity of the DERs in the network is 2180 kW, distributed unevenly across 7 DERs. For each DER, the ramping limit $\overline{o}_g$ in \eqref{eq:subprob} is set at 30\% of its capacity $\overline{s}_g$. Block weighting parameters $\alpha_l$ in~\eqref{eq:master} are set to $10 + 0.01(|\mathcal{D}_l|) \ \forall l \in \mathcal{B}$, where $\mathcal{D}_l$ represents the set of loads in block $l$. The weighting parameter $\omega$ in~\eqref{eq:enumerate} is empirically set to 705 to facilitate convergence.

\subsection{RPOP's performance for varying levels of load uncertainty}
\label{subsec:varying_uncertainty}
\begin{table}
\renewcommand{\arraystretch}{1.0}
 \setlength{\tabcolsep}{4pt}
\centering
\caption{Sensitivity of Operating Costs and Algorithm Run Times to Load Uncertainty Percentage}
\label{uncertainty_sensitivity}
\begin{tabular}{ >{\centering}p{0.10\columnwidth}>{\centering}p{0.10\columnwidth}>{\centering}p{0.10\columnwidth}> {\centering}p{0.11\columnwidth}> {\centering}p{0.10\columnwidth}>
{\centering}p{0.10\columnwidth}>
{\centering\arraybackslash}p{0.12\columnwidth}}
\toprule[1.2pt]\midrule[0.3pt]
Uncertain Load (\%) & Gen Cost (\$) & Load Shed Cost (\$) & \# of Energized Blocks & \# of Closed Switches & Run Time \par (sec.) & Run Time for best IFS (sec.)\\
\midrule
0 & 15,399 & 21,100 & 5 & 4 & 100.8 & 45.67 \\
10 & 12,295 & 31,140 & 4 & 3 & 1952.6 & 1708.2 \\
15 & 12,299 & 31,140 & 4 & 3 & 2500.4 & 2359.8 \\
25 & 9,317 & 41,160 & 3 & 2 & 4599.9 & 4572.5 \\
\midrule[0.3pt]\toprule[1.2pt]
\end{tabular}
\vspace*{-4mm}
\end{table}
Results in this section (from Algorithm~\ref{alg:cutting_plane}) address the 2-stage RPOP at varying load uncertainty levels: 0\% (nominal case), 10\%, 15\%, and 25\%. Table~\ref{uncertainty_sensitivity} details generation and load shedding costs, the number of energized blocks and closed switches, and Algorithm~\ref{alg:cutting_plane} run times. We observe a monotonic increase in the total cost (generation + load shedding  in~\eqref{eq:master}) with higher load uncertainty levels. Moreover, increased uncertainty correlates with reduced energized blocks, lowering generation costs but increasing load shedding costs.

Table~\ref{uncertainty_sensitivity} shows increased run times for cases with load uncertainty compared to the nominal case. This is mainly because obtaining the worst-case uncertainty realization requires solving $2^{\mathcal{|D|}}$ subproblems~\eqref{eq:enumerate} for every iteration of Algorithm~\ref{alg:cutting_plane}. Moreover, run times increase for higher uncertainty levels  requiring more iterations/cuts in Algorithm~\ref{alg:cutting_plane} to approximate the second-stage cost function and 
ensure robustness against larger uncertainty sets. The last column of Table~\ref{uncertainty_sensitivity} displays the run times at which Algorithm~\ref{alg:cutting_plane} reaches the best integer feasible solution (IFS), i.e., the optimal integer solution for the 2-stage RPOP, under each uncertainty level.  Notably, the best IFS is achieved within 81\% of the total run time, on average. 
\begin{figure}[h]
    \centering
    \subfigure[0\% load uncertainty (nominal). Served load = 1190~kW.]{
	\includegraphics[width=3.42in] {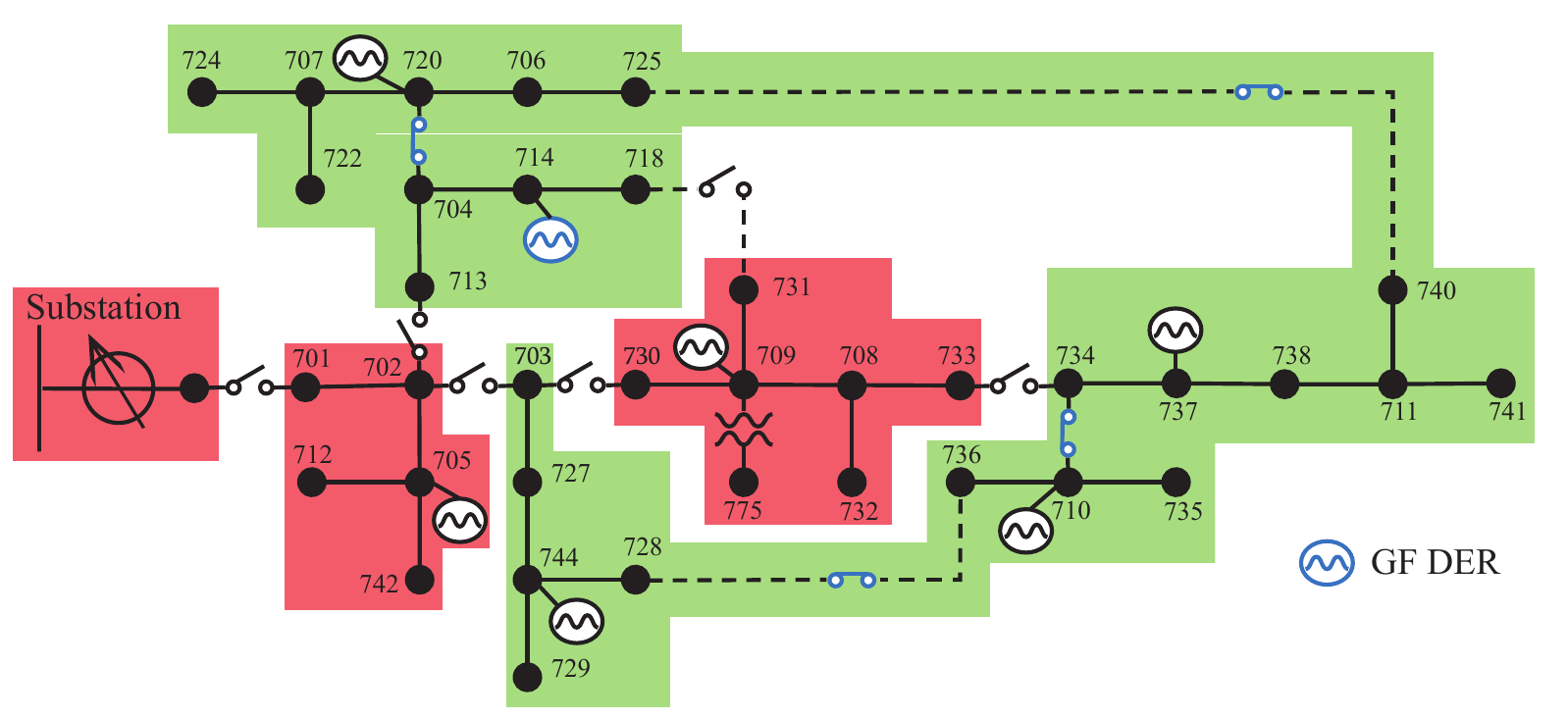}}\\ 
 \vspace{-.2cm}
    \subfigure[10\% and 15\% load uncertainty. Served load =  938 kW.]{
	\includegraphics[width=3.42in]{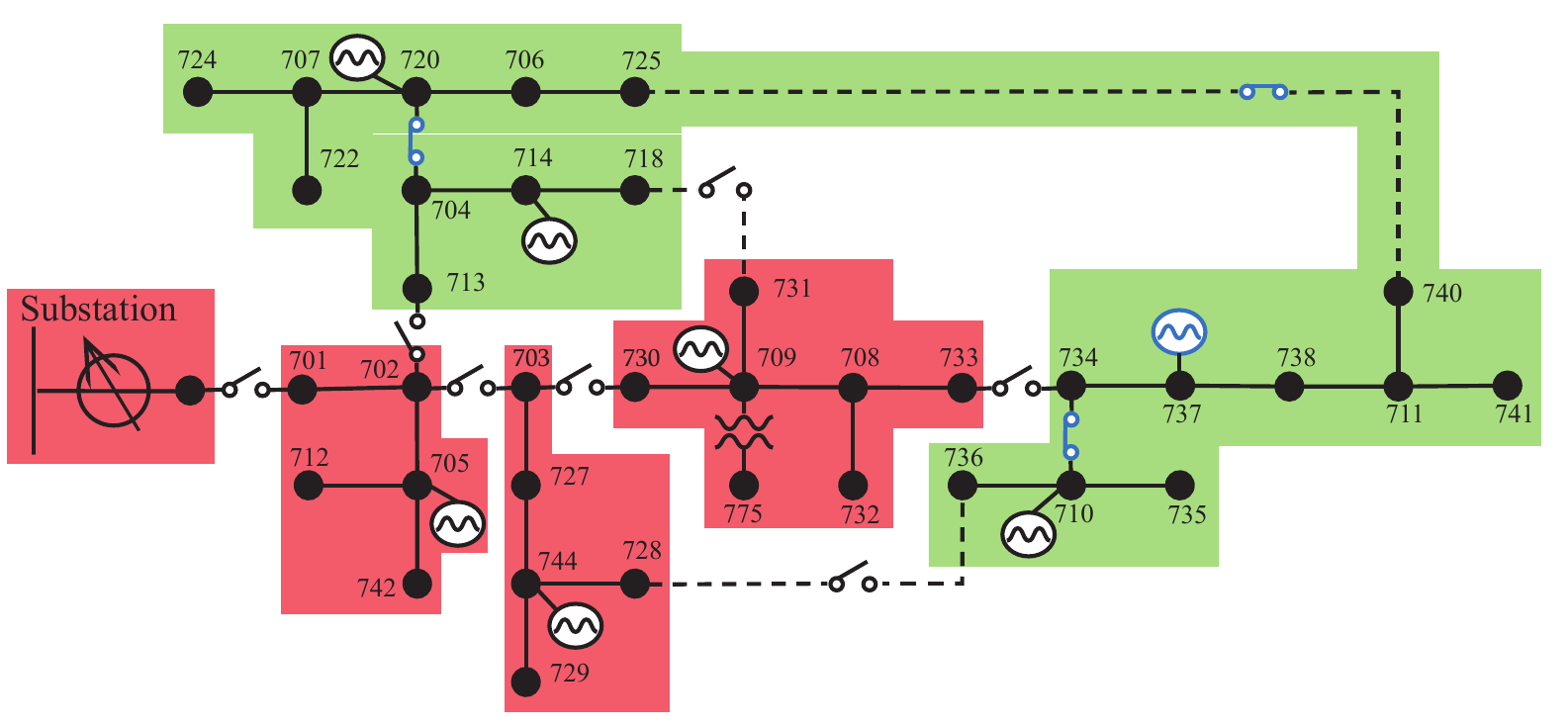}} \\ 
 \vspace{-.2cm}
	\subfigure[25\% load uncertainty. Served load = 811 kW.]{
	\includegraphics[width=3.42in]{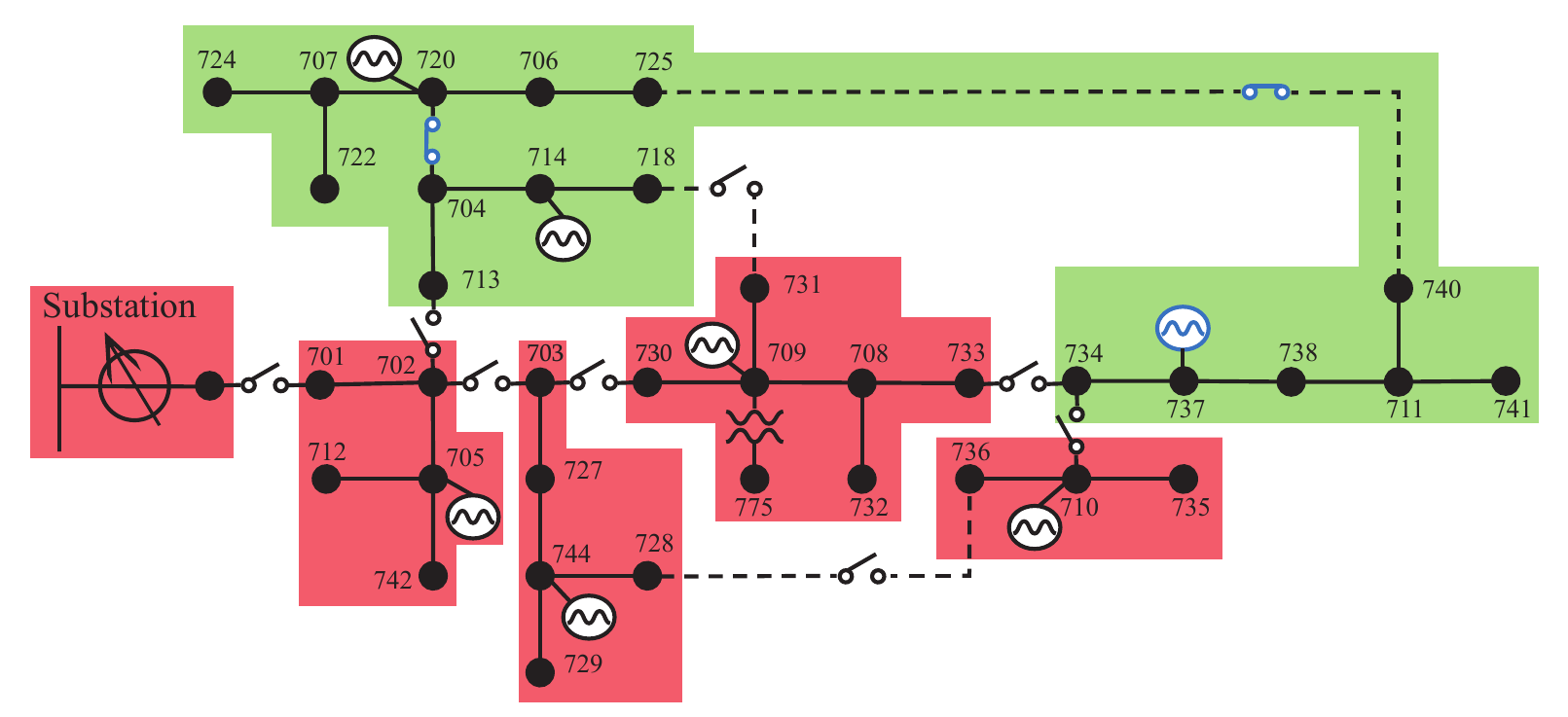}}
    \caption{Optimal partitioning of the 37-bus network with varying levels of uncertain loads. We assume the substation is out-of-service. Energized and de-energized CCs are colored green and red, respectively. GF-DER and closed switches are highlighted in blue.}
    \label{fig:uncertain_37bus_network}
    \vspace*{-4mm}
\end{figure}

Optimal network configurations and block states for various uncertainty levels are shown in Fig.~\ref{fig:uncertain_37bus_network}.
In Fig.~\ref{fig:uncertain_37bus_network}(a), a single inter-connected microgrid is formed using two redundant lines to maximize load delivery. Notably, even without uncertainty,  DERs are unable to satisfy every load.
At 10\% load uncertainty and beyond, satisfying all loads within the energized CC of the nominal case becomes infeasible. To maintain robustness against the worst-case uncertain load, more blocks must be de-energized, as shown in Figs.~\ref{fig:uncertain_37bus_network}(b) and \ref{fig:uncertain_37bus_network}(c).

Without network partitioning (i.e., when the network is fixed to Fig.~\ref{fig:37bus_network}), the optimal solution for all uncertainty levels is to de-energize all loads. This is primarily due to the assumption that individual loads cannot be shed in the distribution network, a practice followed by most electric utilities. In contrast, partitioning the grid and forming \textit{networked microgrids can significantly increase the amount of load served}. For instance, under nominal loading conditions, the optimal partitioning solution of the 2-stage RPOP increases the served load to 1190 kW, representing 53.2\% of the total load.

\subsection{RPOP's performance against deterministic counterparts} 
\label{subsec:deterministic}
We now analyze network partitioning without considering load uncertainty. For this analysis, we fix the network configuration (switch and load block states) to the optimal solution of the nominal load case (as shown in Fig.~\ref{fig:uncertain_37bus_network}(a), and set DER limits as $\pm 30\%$ of their capacity around the established nominal set-points. This limit is assumed reasonable for fast ramping DERs commonly encountered in distribution grids. Subsequently, we assess operational feasibility for varying levels of load uncertainty via out-of-sample feasibility tests. For each level of uncertainty, we generate 10,000 random load samples from the uncertainty set and, for each sample, compute the \textsc{LinDist3Flow} power flow to ascertain feasibility. A sample is deemed feasible if all engineering limits are met.

\textit{The percentage of out-of-sample loads that were found feasible for 10\%, 15\% and 25\% level of load uncertainty were 74.75\%, 67.03\% and 59.96\%}, respectively. As expected, the percent of feasible load samples decreases with increases in uncertainty. 
In summary, this analysis suggests that even with fast ramping DERs to handle real-time load fluctuations, deterministic planning may not always satisfy engineering limits, even at low levels of uncertainty. At higher uncertainty levels, \textit{network partitioning and operation without considering uncertainty can significantly reduce grid resilience} and increase vulnerability to disruptive extreme-weather events.

{\color{black} \subsection{Robust feasibility for non-convex AC power flow}
\label{subsec:acpf_feasibility}
We empirically assess the robustness of the two-stage RPOP's solutions, constrained by linearized \textsc{LinDist3Flow}, against true AC Power Flow (ACPF) constraints, which are non-convex, in three-phase unbalanced grids. Algorithm~\ref{alg:cutting_plane} does not straightforwardly extend for subproblems with non-convex constraints. Instead, for each uncertainty level, we obtain optimal network configurations and DER set-points from Algorithm~\ref{alg:cutting_plane} based on the \textsc{LinDist3Flow} approximation. We fix this network configuration (by fixing all binaries) and set DER capacity limits to  $\pm 30\%$ around the obtained DER set-points. Then, for each uncertainty level, we generate 10,000 random load samples from uncertainty set $\mathcal{U}$, as the worst-case load scenario in the ACPF setting may not occur at an extreme point, i.e., in set $\widehat{\mathcal{U}}$, as discussed in \ref{subsec:2stage_RPOP}. Next, for each level of load uncertainty, we compute the percentage of random load samples for which the ACPF (polar) formulation is feasible. Table~\ref{tab:acpf_robustness} illustrates that very few load samples render the ACPF formulation infeasible, even at higher uncertainty levels. To summarize, for the IEEE 37-bus case, the \textsc{LinDist3Flow} formulation in the two-stage RPOP provides a sufficiently robust approximation, being almost always feasible under the non-convex ACPF constraints.}
\begin{table}
\renewcommand{\arraystretch}{1.0}
\centering
\caption{\color{black} Robust feasibility of two-stage RPOP solutions against non-convex three-phase AC power flow at varying levels of uncertainty.}
\label{tab:acpf_robustness}
\begin{tabular}{ >{\centering}p{0.4\columnwidth}>
{\centering\arraybackslash}p{0.4\columnwidth}}
\toprule[1.2pt]\midrule[0.3pt]
Load Uncertainty (\%) &  AC feasible load samples (\%) \\
\midrule
10 & 100\\
15 & 98.23\\
25 & 99.69 \\
\midrule[0.3pt]\toprule[1.2pt]
\end{tabular}
\vspace*{-3mm}
\end{table}

\section{Conclusions}
In this paper, we introduced a mixed-integer robust partitioning and operation problem (RPOP) to optimize network configuration and generator set-points, minimizing generation costs and undelivered uncertain load. To make this problem tractable, we proposed a novel two-stage RPOP reformulation. The first stage optimizes network partitioning decisions using a mixed-integer program to maximize load delivery, while the second stage verifies robustness and feasibility of linearized three-phase unbalanced power flow. Additionally, we presented a cutting-plane algorithm for efficient solving of this two-stage RPOP. Illustrating the algorithm's efficacy during a disruptive event, we conducted a case study using a modified IEEE 37-bus test system. Results demonstrated that partitioning a distribution grid into networked microgrids significantly reduces load shedding. Furthermore, the optimal planning decisions under nominal loading conditions were found to become infeasible to power flow under worst-case load uncertainty, even at uncertainty levels as low as 10\%. This underscores the necessity of a robust formulation. Finally, despite relying on a \textsc{LinDist3Flow} approximation for robustness, the method's planning and operational solutions remained adequately accurate in the non-convex AC framework.

In future work, we aim to scale these approaches for larger and realistic distribution grids, utilizing load-clustering-based heuristics to expedite solving the $\max-\min$ problems. Additionally, we will explore less conservative uncertainty sets and integrate battery storage devices and renewable sources.

\bibliographystyle{IEEEtran}
\bibliography{references.bib}

\begin{thebibliography}{10}
\providecommand{\url}[1]{#1}
\csname url@samestyle\endcsname
\providecommand{\newblock}{\relax}
\providecommand{\bibinfo}[2]{#2}
\providecommand{\BIBentrySTDinterwordspacing}{\spaceskip=0pt\relax}
\providecommand{\BIBentryALTinterwordstretchfactor}{4}
\providecommand{\BIBentryALTinterwordspacing}{\spaceskip=\fontdimen2\font plus
\BIBentryALTinterwordstretchfactor\fontdimen3\font minus
  \fontdimen4\font\relax}
\providecommand{\BIBforeignlanguage}[2]{{%
\expandafter\ifx\csname l@#1\endcsname\relax
\typeout{** WARNING: IEEEtran.bst: No hyphenation pattern has been}%
\typeout{** loaded for the language `#1'. Using the pattern for}%
\typeout{** the default language instead.}%
\else
\language=\csname l@#1\endcsname
\fi
#2}}
\providecommand{\BIBdecl}{\relax}
\BIBdecl

\bibitem{Hamidieh2022Microgrids}
M.~Hamidieh and M.~Ghassemi, ``Microgrids and resilience: A review,''
  \emph{IEEE Access}, vol.~10, pp. 106\,059--106\,080, 2022.

\bibitem{Li2017Networked}
Z.~Li, M.~Shahidehpour, F.~Aminifar, A.~Alabdulwahab, and Y.~Al-Turki,
  ``Networked microgrids for enhancing the power system resilience,''
  \emph{Proceedings of the IEEE}, vol. 105, no.~7, pp. 1289--1310, 2017.

\bibitem{barnes2019resilient}
A.~Barnes, H.~Nagarajan, E.~Yamangil, R.~Bent, and S.~Backhaus, ``Resilient
  design of large-scale distribution feeders with networked microgrids,''
  \emph{Electr Pwr Syst Res}, vol. 171, pp. 150--157, 2019.

\bibitem{Golari2014Two}
M.~Golari, N.~Fan, and J.~Wang, ``Two-stage stochastic optimal islanding
  operations under severe multiple contingencies in power grids,''
  \emph{Electric Power Systems Research}, vol. 114, pp. 68--77, 2014.

\bibitem{fobes2022optimal}
D.~M. Fobes, H.~Nagarajan, and R.~Bent, ``Optimal microgrid networking for
  maximal load delivery in phase unbalanced distribution grids: A declarative
  modeling approach,'' \emph{IEEE Transactions on Smart Grid}, vol.~14, no.~3,
  pp. 1682--1691, 2022.

\bibitem{Zhou2022Three}
A.~Zhou, H.~Zhai, M.~Yang, and Y.~Lin, ``Three-phase unbalanced distribution
  network dynamic reconfiguration: A distributionally robust approach,''
  \emph{IEEE Transactions on Smart Grid}, vol.~13, no.~3, pp. 2063--2074, 2022.

\bibitem{Lee2015Robust}
C.~Lee, C.~Liu, S.~Mehrotra, and Z.~Bie, ``Robust distribution network
  reconfiguration,'' \emph{{IEEE} Transactions on Smart Grid}, vol.~6, no.~2,
  pp. 836--842, 2015.

\bibitem{Mahdavi2023Robust}
M.~Mahdavi, K.~E.~K. Schmitt, and F.~Jurado, ``Robust distribution network
  reconfiguration in the presence of distributed generation under uncertainty
  in demand and load variations,'' \emph{IEEE Transactions on Power Delivery},
  vol.~38, no.~5, pp. 3480--3495, 2023.

\bibitem{babaei2020distributionally}
S.~Babaei, R.~Jiang, and C.~Zhao, ``Distributionally robust distribution
  network configuration under random contingency,'' \emph{IEEE Transactions on
  Power Systems}, vol.~35, no.~5, pp. 3332--3341, 2020.

\bibitem{barani2018optimal}
M.~Barani, J.~Aghaei, M.~A. Akbari, T.~Niknam, H.~Farahmand, and M.~Korp{\aa}s,
  ``Optimal partitioning of smart distribution systems into supply-sufficient
  microgrids,'' \emph{IEEE Transactions on Smart Grid}, vol.~10, no.~3, pp.
  2523--2533, 2018.

\bibitem{arefifar2012supply}
S.~A. Arefifar, Y.~A.-R.~I. Mohamed, and T.~H. El-Fouly,
  ``Supply-adequacy-based optimal construction of microgrids in smart
  distribution systems,'' \emph{IEEE Trans. on smart grid}, vol.~3, no.~3, pp.
  1491--1502, 2012.

\bibitem{Gholami2019Proactive}
A.~Gholami, T.~Shekari, and S.~Grijalva, ``Proactive management of microgrids
  for resiliency enhancement: An adaptive robust approach,'' \emph{IEEE Trans.
  on Sustainable Energy}, vol.~10, no.~1, pp. 470--480, 2019.

\bibitem{robbins2015optimal}
B.~A. Robbins and A.~D. Dom{\'\i}nguez-Garc{\'\i}a, ``Optimal reactive power
  dispatch for voltage regulation in unbalanced distribution systems,''
  \emph{IEEE Transactions on Power Systems}, vol.~31, no.~4, pp. 2903--2913,
  2015.

\bibitem{Claeys2021No}
S.~Claeys, F.~Geth, M.~Sankur, and G.~Deconinck, ``No-load linearization of the
  lifted multi-phase branch flow model: Equivalence and case studies,'' in
  \emph{2021 {IEEE} {PES} Innovative Smart Grid Technologies Europe ({ISGT}
  Europe)}, 2021, pp. 1--5.

\bibitem{yang2021robust}
H.~Yang, D.~P. Morton, C.~Bandi, and K.~Dvijotham, ``Robust optimization for
  electricity generation,'' \emph{INFORMS Journal on Computing}, vol.~33,
  no.~1, pp. 336--351, 2021.

\bibitem{Gan2014Convex}
L.~Gan and S.~H. Low, ``Convex relaxations and linear approximation for optimal
  power flow in multiphase radial networks,'' in \emph{Power Systems
  Computation Conference}, 2014, pp. 1--9.

\bibitem{Fobes2020PMD}
D.~M. Fobes, S.~Claeys, F.~Geth, and C.~Coffrin,
  ``{PowerModelsDistribution.jl}: An open-source framework for exploring
  distribution power flow formulations,'' \emph{Electric Power System
  Research}, vol. 189, p. 106664, 2020.

\bibitem{Lei2020Radiality}
S.~Lei, C.~Chen, Y.~Song, and Y.~Hou, ``Radiality constraints for resilient
  reconfiguration of distribution systems: Formulation and application to
  microgrid formation,'' \emph{IEEE Transactions on Smart Grid}, vol.~11,
  no.~5, pp. 3944--3956, 2020.

\bibitem{Li2022Revisiting}
Y.~Li, Y.~Gu, and T.~C. Green, ``Revisiting grid-forming and grid-following
  inverters: A duality theory,'' \emph{IEEE Transactions on Power Systems},
  vol.~37, no.~6, pp. 4541--4554, 2022.

\bibitem{Sadeque2021Multiple}
F.~Sadeque, D.~Sharma, and B.~Mirafzal, ``Multiple grid-forming inverters in
  black-start: The challenges,'' in \emph{IEEE 22nd Workshop on Control and
  Modelling of Power Electronics (COMPEL)}, 2021, pp. 1--6.

\bibitem{Han2016Review}
H.~Han, X.~Hou, J.~Yang, J.~Wu, M.~Su, and J.~M. Guerrero, ``Review of power
  sharing control strategies for islanding operation of ac microgrids,''
  \emph{IEEE Trans. on Smart Grid}, vol.~7, no.~1, pp. 200--215, 2016.

\bibitem{Geoffrion1972Generalized}
A.~Geoffrion, ``Generalized benders decomposition,'' \emph{Journal of
  Optimization Theory and Applications}, vol.~10, no.~4, p. 237–260, 1972.

\bibitem{Bertsekas2009Convex}
D.~P. Bertsekas, \emph{Convex Optimization Theory}.\hskip 1em plus 0.5em minus
  0.4em\relax Athena Scientific, 2009.

\bibitem{dunning2017jump}
I.~Dunning, J.~Huchette, and M.~Lubin, ``Ju{MP}: A modeling language for
  mathematical optimization,'' \emph{SIAM Review}, vol.~59, no.~2, pp.
  295--320, 2017.

\bibitem{gurobi}
\BIBentryALTinterwordspacing
{Gurobi Optimization, LLC}, ``{Gurobi Optimizer Reference Manual},'' 2022.
  [Online]. Available: \url{https://www.gurobi.com}
\BIBentrySTDinterwordspacing

\bibitem{kersting2001radial}
W.~H. Kersting, ``Radial distribution test feeders,'' in \emph{2001 IEEE Power
  Engineering Society Winter Meeting. Conference Proceedings (Cat. No.
  01CH37194)}, vol.~2.\hskip 1em plus 0.5em minus 0.4em\relax IEEE, 2001, pp.
  908--912.

\end{thebibliography}

\end{document}